\documentclass[12pt]{amsart}

\def\checkbox{\leavevmode\vbox to 9pt{\hrule \vss
	\hbox to 9pt{\vrule height 9pt \hfil\vrule height 9pt}\vss
	\hrule}\ }

\newcommand{\Q}{\mathbb Q}
\newcommand{\R}{\mathbb R}
\newcommand{\Z}{\mathbb Z}

\renewcommand{\epsilon}{\varepsilon}
\renewcommand{\phi}{\varphi}
\newtheorem{Lemma}{Lemma}[section]
\newtheorem{Theorem}{Theorem}[section]
\newtheorem{Proposition}{Proposition}[section]

\newtheorem{Corollary}[Proposition]{Corollary}
\newtheorem{Definition}{Definition}[section]

\newtheorem{Remark}{Remark}[section]

\begin{document}
\address{Department of Mathematics, Pennsylvania State University, University Park, PA $16802$, USA
}
\author[]{Alexandr Borisov}
\title[A. Borisov, Positive Positive-definite Functions]{Positive positive-definite functions and measures on locally compact abelian groups}
\email{borisov@math.psu.edu}
\date{\today}
\maketitle
\thispagestyle{empty}
\centerline{Preliminary version}

\section{Introduction}
In the paper \cite{Bor} we gave a cohomological interpretation of Tate's Riemann-Roch formula using some new harmonic analysis objects, ghost-spaces. When trying to investigate these objects in general, we realized the importance of functions and measures on locally compact abelian groups that are both positive and positive-definite at the same time. It looks like this class of functions and measures was not systematically studied before. The goal of this paper is to partially fill in this gap. We answer some of the natural questions involving these functions and measures, especially those that satisfy some extra integrability conditions. We also study some operations and constructions involving these functions and measures.

There are several very interesting open questions, that we are only able to point out at this moment. In particular, the structure of the cone of such functions is not clear even when the group is just $\R.$ Please refer to section 5 where this and other open problems are discussed.

The paper is organized as follows. In section 2 the positive positive-definite functions and measures are defined, some of their properties are discussed, and some simple constructions involving them are carried out. In sections 3 and 4 we restrict our attention to such functions and measures that satisfy some extra integrability conditions. Finally, in section 5 we point out some natural open questions that we were unable to answer.

{\bf Acknowledgments.} The author thanks Jeff Lagarias, Barry Mazur and Christopher Deninger whose interest in \cite{Bor} motivated the author to continue his work in this direction. The author also thanks Joaquim Ortega-Cerda for the references to Hardy's theorem and its generalizations.

\section{General definitions and results}
In this section we are going to define some classes of functions and measures on locally compact abelian groups. We  will prove that these classes correspond to each other via duality and are stable under the operations of pull-back of functions and push-forward of measures, with respect to continuous homomorphisms of groups. They are also stable under addition and multiplication of functions, and convolution of measures.

We are going to use the book of Folland \cite{Fol} as our basic harmonic analysis  reference. In particular, we are going to use his terminology regarding positive-definite functions and functions and measures of positive type. So a positive-definite function is not necessarily continuous. It is just any function $f$ on a locally compact group $G,$ such that for all $n-$tuples $(x_i)^n_{i=1}$ of elements of $G$ and complex numbers $(c_i)^n_{i=1}$
$$\sum \limits_{i,j=1}^{n} c_i \bar{c_j} f(x_i-x_j) \ge 0$$
 However most of our functions are going to be continuous anyway, so this is not important.

Let us recall a little bit of Pontryagin duality theory. A good exposition of it can be found, e.g. in \cite{Fol}. To any locally compact abelian group $G$ one can associate its group of characters, or dual group $\hat{G}.$ The Pontryagin duality theorem tells that $\widehat{\widehat{G}}=G.$ To every bounded Radon measure $\mu$ on $\widehat{G}$ one can associate a function $\check{\mu}$ on $G$ (its inverse Fourier transform) by the following formula.
$$\check{\mu} (x)= \int \limits_{\chi \in \hat{G}} \chi (x) d\mu (\chi )$$
To every $L^1$ function $f$ on $G$ one can associate a function $\hat{f}$ on $\hat{G}$ (its Fourier transform) that actually depends on the Haar measure $m$ on $G$ by the following formula.
$$\hat{f}(\chi)=\int \limits_{x\in G} \overline{\chi(x)}\cdot f(x) dm(x) $$

So here are our basic definitions.

\begin{Definition} 
Suppose $G$ is a locally compact abelian group, and $f$ is a function on $G$. We say that $f$ is a positive positive-definite function (PPD, for short) on $G$  if the following properties are satisfied.

1) $f(x)$ is real-valued and non-negative for all $x\in G$

2) $f$ is of positive type on $G$ (i.e. continuous and positive-definite)
\end{Definition}

 If $f$ is a PPD function, then it is even, i.e. $f(-x)=f(x)$ for all $x$. This follows from $f$ being real-valued and of positive type (cf. \cite{Fol}, prop. 3.22).

\begin{Lemma}
Suppose $G$ is a locally compact abelian group, and $f$ is an PPD function on $G$. Then there is a unique bounded Radon measure $\mu$ on $\widehat{G}$, such that $f= \check{\mu}.$ This measure $\mu$ is nonnegative real-valued measure, which is even and positive-definite.
\end{Lemma}

{\bf Proof.} The existence of $\mu$ follows from the Bochner theorem (cf. \cite{Fol}, prop. 4.18). The uniqueness is easy. The properties of $\mu $ follow easily from the properties of $f.$  \hfill \checkbox

\begin{Definition}
We will say that measure $\mu $ on $G$ is positive positive-definite (PPD for short) if $\check{\mu}$ is a PPD function on $\widehat{G}$. We will call $\check{\mu}$ the function dual to the measure $\mu$. We will call $\mu$ the measure dual to the function $\check{\mu}.$ \end{Definition}

Obviously, every PPD measure has unique PPD function dual to it. By the above lemma every PPD function has unique PPD measure dual to it.

For a locally compact abelian group $G$ we will denote by $PPDf(G)$ the set of PPD functions  and by $PPDm(G)$ the set of PPD measures on it.

The following theorem is proved in \cite{Bor} . We reproduce its proof here for the convenience of the reader.

\begin{Theorem} Suppose $G$ is a locally compact abelian group and $u\in PPDf(G)$. Then for all $x\in G$ $u(x)\le u(0)$. Also, those $x$ that $u(x)=u(0)$ form a closed subgroup $H$ of $G.$  Moreover, $u(x)$ is a pull-back of a PPD function on $G/H$ via the natural morphism $G\rightarrow G/H.$
\end{Theorem}

{\bf Proof.} The first claim is contained in Folland \cite{Fol}, cor. 3.32. To prove the second and third claims we note that by \cite{Fol}, prop. 3.35 the following matrix is positive definite.
$$ \left [\begin {array}{lcr} \mbox {u(0)}&\mbox {u(x)}&\mbox {u(x+y)}
\\\noalign{\medskip}\mbox {u(x)}&\mbox {u(0)}&\mbox {u(y)}
\\\noalign{\medskip}\mbox {u(x+y)}&\mbox {u(y)}&\mbox {u(0)}\end {array}
\right ]
 $$
If $u(x)=u(0)$, it implies that $(u(x+y)-u(y))^2\le 0,$ so $u(x+y)=u(y).$ This implies that $u$ is a pull-back of some function on $G/H.$ It is now a trivial check to establish that this function on $G/H$ is PPD. \hfill \checkbox

If $f$ is a PPD function on $G$ then so is $\alpha \cdot f$ for any  positive real $\alpha .$ A lot of our constructions involve Fourier transform which depends on the choice of a Haar measure. So the results of such operations are only defined up to a constant. For this reason, it is convenient to consider "normalized" functions and measures as in the definitions below. We should also point out that this normalization condition naturally appears in \cite{Bor} in the definitions of ghost-spaces.

\begin{Definition} A PPD function u on $G$ is normalized if $u(0) =1.$ We will denote the set of all normalized PPD functions on $G$ by NPPDf(G). When we say that we normalize a PPD function, this means that we multiply it by a suitable constant to make it normalized.
\end{Definition}

\begin{Definition} A PPD measure $\mu $ on $G$ is normalized if $\mu$ is a probability measure. We will denote the set of all normalized PPD measures on $G$ by NPPDm(G). When we say that we normalize a PPD measure, this means that we multiply it by a suitable constant to make it normalized.
\end{Definition}

\begin{Remark} A measure $\mu \in PPDm(\widehat{G})$ is normalized if and only if the function $\check{\mu}\in PPDf(G)$ is normalized.
\end{Remark}

\begin{Remark}
Suppose $\pi: H\rightarrow G$ is a continuous homomorphism of locally compact abelian groups.
Denote by $\pi^*$ the pull-back of functions and by $\pi_*$ the push-forward of measures.
Then the following is always true.

1) If $f\in NPPDf(G)$ then $\pi^*f \in NPPDf(H)$.

2) If $\mu \in NPPDm(H)$ then $\pi_* \mu \in NPPDm(G)$.
\end{Remark}

Here are some examples of PPD functions and measures. More examples can be obtained from these using pull-backs and push-forwards. Also the set of PPD functions (or measures) on a given group is obviously a convex cone. It is also closed under multiplication as will be shown in Corollary 4.2. Finally, for some PPD functions and measures (which we will later call good, cf. Definitions 3.1 and 3.2) there are other constructions available, like a corestriction of functions, restriction of measures, etc. Please refer to section 4 for the details.

{\bf Examples.}

1) The constant $1$ is a PPD function for every $G$.

2) The point measure at $0$ is a PPD measure for every $G$.

3) For every positive-definite quadratic form $Q$ on $\R^n$ (or $\Z^n$) the function $e^{-Q(x,x)}$ is a PPD function on $\R^n$(or $\Z^n$). And if we multiply it by a Haar measure on the corresponding group, we get a PPD measure.

\section{Good functions and measures}
The class of (normalized) PPD functions and measures was used in \cite{Bor} for the definitions of ghost-spaces. However it turns out that in order to have good categorical properties of ghost-spaces one has to restrict the class of functions and measures used. This was the original  reason for introducing the class of good functions and measures. It turned out afterwards that this class of functions and measures is pretty natural. It is stable under some interesting constructions and is much more interesting than the class of all PPD functions and measures. Please refer to section 4 for the details.

\begin{Definition} 
Suppose $G$ is a locally compact abelian group, and $f$ is a function on $G$. We say that $f$ is good if the following properties are all satisfied.

1) $f\in PPDf(G).$

2) $f \in L^1(G)$ with respect to a Haar measure on $G.$

3) The Fourier transform $\hat{f}$ of $f$ is continuous and belongs to $L^1(\widehat{G}).$

4) Both $f$ and $\hat{f}$ are strictly positive.

5) Both $f$ and $\hat{f}$ are $L^1$ with respect to a Haar measure when restricted
to arbitrary closed subgroups of $G$ or $\widehat{G}$ respectively.
\end{Definition}

\begin{Remark}

1) By the Fourier Inversion theorem $f=\Check{(\hat{f} m)},$ where $m$ is a suitably normalized Haar measure on $\widehat{G}$ (cf., e.g. \cite{Fol}).

2) If $f$ is good then so is $\hat{f}$.

3) The condition 5 doesn't follow from the others, as the following example shows.

\end{Remark}

{\bf Example.} Consider $G=\R ^2,$ with coordinates $(x,y).$
Consider the following function $f$ on $G.$
$$f(x,y) = \sum \limits_{n=1}^{\infty} \frac{1}{n^2} e^{-\pi (n^2 x^2 + \frac{1}{n^2} y^2)}$$

This sum converges for all $(x,y)$ to a continuous function,
which is $L^1$ on $G$. Its Fourier transform is up to a constant $f(y,x)$,
so it has the same properties. However the restriction of $f$ to the line $x=0$ is not in $L^1.$

Some interesting properties of good functions are summarized in the following lemma.

\begin{Lemma}
Suppose $G$ is a locally compact abelian group, and $f$ is a good function on it. Then

1) $0 < f(x) \le f(0)$ for all $x\in G$

2) $f(x)=f(0)$ implies $x=0$

\end{Lemma}

{\bf Proof.} The first statement is proven in Theorem 2.1 for the more general class of PPD functions. It is also proven there that the set of all elements $x$ of $G$ such that $f(x)=f(0)$ is a closed subgroup of $G.$ Moreover it is the stabilizer of $f$ in $G$, i.e. it consists of the elements $x\in G$ such that $f$ is invariant under translation by $x$. If this subgroup, say $H,$ is not trivial, then  $\hat{f}(\chi) =0$ for all $\chi$ that are not trivial on $H.$ This would contradict the condition 4 of the Definition 3.1. \hfill \checkbox

\begin{Definition}
Suppose $G$ is a locally compact abelian group, and $\mu$ is a Radon measure on $G$. We say that $\mu$ is good if $\mu =f\cdot m,$ where $f$ is a good function and $m$ is a Haar measure on $G.$
\end{Definition}

 For a locally compact abelian group $G$ we will denote by $Goodf(G)$ the set of good functions  and by $Goodm(G)$ the set of good measures on it.

\begin{Lemma} Suppose $G$ is a locally compact abelian group.

1) If $\mu$ is a good measure on $\widehat{G}$ then $\check{\mu} $ is a good function on $G$.

2) If $f$ is a good function on $G$ then there exists a unique measure $\mu$ on $\widehat{G}$ such that $f=\check{\mu}$ and this measure $\mu$ is good.
\end{Lemma}

{\bf Proof.} The proofs are pretty straightforward. The measure $\mu$ in the second statement can be obtained by multiplying $\hat{f}$ by a suitably normalized Haar measure. \hfill \checkbox

Naturally, we will denote by $NGoodf(G)$ the set of all normalized good functions on $G$, and by $NGoodm(G)$ the set of all normalized good measures on $G$, cf. Definitions 2.3 and 2.4. 

\begin{Definition} If $f\in NGoodf(G)$ then its {\bf normalized dual} is a function $\hat{f} \in NGoodf(\widehat{G})$ which is equal to the Fourier transform of $f$ for some Haar measure on $G.$
\end{Definition}

To any normalized good function $f$ one can associate the unique Haar measure $m$ on $G$ such that $f\cdot m$ is a normalized measure. Then the Fourier transform of $f$ with respect to this measure is a normalized dual of $f$. Dually, to any normalized good measure $\mu$ on $\widehat{G}$ one can associate the unique Haar measure $\hat{m}$ on $\widehat{G}$ such that $\mu = \hat{f}\cdot \hat{m},$ where $\hat{f} \in NPPDf(\widehat{G}).$ If $f=\check{\mu}$ then $m$ and $\hat{m}$ are dual Haar measures. This easy fact is behind Theorem 4.1 of \cite{Bor}.

\section{Some constructions involving good functions and measures}
In this section we define some interesting constructions involving good functions and measures.
Because good measures are dual to good functions, we will do most of our constructions for the functions. And we will leave it to the reader to define similar operations for the measures via duality. 

\begin{Definition}
Suppose $H$ is a closed subgroup of a locally compact abelian group $G$ and $f\in NGoodf(G).$ Then we call the pull-back of $f$ via the natural embedding $H\rightarrow G$ the {\bf restriction} of $f$ to $H$. We denote it by $f_{|_H}$.
\end{Definition}

\begin{Definition} Suppose $K=G/H$ is the quotient of a locally compact abelian group $G$ by a closed subgroup $H$. Suppose $f \in NGoodf(G).$ Then we define {\bf corestriction} of $f$ to $K$ as follows. We take the Fourier dual $\hat{f}$ of $f$ and restrict it to $H^{\perp}=\widehat{K}$. Then we take its Fourier dual normalized so that its value at $0$ is $1.$ This function is what we call the corestriction of $f$ to $K$, to be denoted by $f^{|^K}.$
\end{Definition}

The following theorem is our main result about restrictions and corestrictions of good functions.

\begin{Theorem} Suppose $H$ is a closed subgroup of a locally compact abelian group $G$ and $f \in Goodf(G).$ Then the following is always true.

1) The restriction $f_{|_H} \in NGoodf (H)$

2) The corestriction of $f^{|^{G/H}} \in NGoodf(G/H)$

3) If $g=f^{|^{G/H}}$ is the above corestriction, and $m_H$ is a Haar measure on $H,$ then for almost all $x\in G$
$$g(\bar{x})=\frac{\int \limits_{y \in H}f(x+y) dm_h(y)}{\int \limits_{y\in H}f(y) dm_H (y)}$$
Also for $x=0$ the above equality is always true, and for {\bf all} $x$ the left hand side is not less than the right hand side.

\end{Theorem}

{\bf Proof.} Let us first choose some Haar measures $m_H$ and $m_G$ on $H$ and $G.$ Then let us denote by $m_{\widehat{G}}$ the dual measure to $m_G$ and by
$m_{\widehat{G}/H^{\perp}}$ the dual measure to $m_H$ ). Then we denote by $m_{H^{\perp}}$ such Haar measure on $H^{\perp}$ that the Fubini theorem is satisfied for $m_{H^{\perp}}$, $m_{\widehat{G}/H^{\perp}}$, and $m_{\widehat{G}}$. Then if $m_{G/H}$ is the dual measure to $m_{H^{\perp}}$,the measures $m_H$, $m_{G/H}$, and  $m_G$ satisfy the Fubini theorem (cf. \cite{Fol}).

 The main part of the proof of the theorem is the following proposition.

\begin{Proposition} Under the above notations and assumptions, let us define two
functions on $G/H$.

$$u(\bar{x}) = \int \limits_{\chi \in H^{\perp}} \hat{f} (\chi) \chi (x) dm_{H^{\perp}} (\chi)$$

$$v(\bar{x}) = \int \limits_{y \in H} f(x+y) dm_H (y)$$

Then 

1) $v(\bar{x})$ converges for all $x.$

2) $v(\bar{x}) \le u(\bar{x})$ for all $x.$

3) $v(0) =u(0).$
\end{Proposition}

{\bf Proof.} First of all, it is easy to see that $v(\bar{x})$ converges to the continuous function $u(\bar{x})$ everywhere except possibly for a set of Haar measure zero. To prove the first statement suppose for some $x_0$ $v(\bar{x_0})$ doesn't converge. Because $f$ is positive and continuous this implies that for some big enough open set with compact closure $U\subset H$
$$\int \limits_{y\in U} f(x_0+y) dm_H(y) > A,$$
where $A=u(0).$ 

Because the set of the points where $v(\bar{x})\neq u(\bar{x})$ has measure zero, there is a sequence of points ${x}$ converging to $x_0$ such that $v(\bar{x})=u(\bar{x})$.
Then by the continuity of $f$ for $x$ close enough to  $x_0$ 
$$u(\bar{x}) \ge \int \limits_{y \in U} f(x+y) dm_H(y) > A$$
But this is impossible by Theorem 2.1.

The same argument obviously implies the second statement. To prove the third statement, consider $A$ as above and also $B=v(0).$ We just proved that $A \ge B.$ But 
$$A=\int \limits_{\chi \in H^{\perp}} \hat{f} (\chi ) dm_{H^{\perp}}(\chi )$$
and
$$B=\int \limits_{y \in H} f(y) dm_H(y)$$
So we can apply the above argument dually and get that $B\ge A.$ This proves the proposition.

The  theorem  now follows easily. We leave the details to the reader. \hfill \checkbox

One can similarly define corestriction and restriction of normalized good measures. Namely, the corestriction to the quotient is just the push-forward. And the restriction to a subgroup can be defined by restricting the corresponding good function and then multiplying by a suitable Haar measure of a subgroup. One can easily check that these definitions are transitive and compatible with Pontryagin duality.

\begin{Theorem}
Suppose $u\in PPDf(G),$ and $v\in PPDf(H)$. Consider the function $w$ on $G\bigoplus H$ defined as below. Here $\pi_G$ and $\pi_H$ are the projections from $G \bigoplus H$ onto $G$ and $H.$
$$w=(\pi_G)^*(u) \cdot (\pi_H)^*(v)$$
Then $w\in PPDf(G\bigoplus H).$ Also if $u$ and $v$ are actually good than $w$ is also good. And if $u$ and $v$ are normalized then $w$ is also normalized.
\end{Theorem}

{\bf Proof.} To prove that $w$ is PPD, we need to check all the conditions in the Definition 2.1. The hardest one is that $w$ is positive-definite. This follows from the fact that $w$ is the inverse Fourier transform of the convolution of two measures on $\widehat{G} \bigoplus \widehat{H}$ that are push-forwards of the measures dual to $u$ and $v$.

So we assume that $u$ and $v$ are good and we need to prove that $w$ is also good. It is pretty easy to check all conditions from Definition 3.1 except for the condition 5. For the condition 5, we should check the $L^1$ property for the restrictions of $f$ and restrictions of $\hat{f}.$ We will prove this for $f,$ the statement for $\hat{f}$ then follows by duality.

So, suppose $M$ is a closed subgroup of $G\bigoplus H.$ We need to prove that the restriction of $f$ to $M$ is $L^1$ with respect to a Haar measure on $M.$ By Theorem 4.1 (1), we can reduce the problem to the case when $\pi_G(M)=G$ and $\pi_H(M)=H.$ So this is now our assumption.
If $N= M \cap G $ and $K= \pi_H (\pi_G^{-1}(N))$ we can define a map $\rho : G \rightarrow (H/K)$ as follows.
$$\rho(x) = \pi_H (y) \cdot K,$$
where $\pi_G(y) = x$.
Then
$$\int \limits_{(x,y)\in M} w(x,y)dm_M (x,y)= const \cdot \int \limits_{x\in G} u(x)\cdot f(\rho (x)) dm_G(x),$$
where $f$ is a function on $H/K$ such that for $y\in H$

$$f(y)=\int \limits_{z\in K}  v(y+z) dm_K(z) $$
By Theorem 4.1 (3) $f \le f(0) \cdot v^{|^{H/K}}$. Therefore  $f(\rho(x)) \le f(0)$ for all $x\in G.$ Thus the  above integral is bounded by
$$f(0)\cdot \int \limits_{x\in G} u(x) dm_G(x) < \infty$$
\hfill \checkbox

\begin{Corollary} The product of two PPD (good) functions on a group $G$ is PPD (good). The convolution of two PPD (good)  measures on $G$ is PPD (good). Also if the functions (measures) are normalized then so is their product (convolution).
\end{Corollary}

{\bf Proof.} If $u$ and $v$ are two functions on $G$ then   $u\cdot v$ is just the restriction to the diagonal of the function $w$ on $G\bigoplus G$. So the statements for the functions follow from the above theorem. The statements for the measures follows  via duality.
\hfill \checkbox

In fact,  stronger results can be obtained. Namely, the product of two PPD functions or convolution of two PPD measures is good as soon as at least one of the original function or measures is good. The following theorem proves it for the functions. The statement for the measures follows via duality.

\begin{Theorem} Suppose $G$ is a locally compact abelian group, $f\in NPPDf(G),$ and $g\in Goodf(G).$ Then $w=f\cdot g \in Goodf(G).$
\end{Theorem}

{\bf Proof.} By Corollary 4.2 we only need to check some integrability conditions for $w$ and $\hat{w}.$  By Theorem 2.1 $w(x)\le g(x)$ for all $x\in G.$ So the integrability conditions for $w$ follow from the corresponding integrability conditions for $g.$ We just need to prove that $\hat{w}$ is $L^1$ when restricted to any closed subgroup $H$ of $\widehat{G}.$ Suppose $\mu \in PPDm(\widehat{G})$ is a measure dual to $f$, and $u\in NGoodf(\widehat{G})$ is the normalized dual of $g.$ Then the function $\widehat{w}$ is proportional to the function $u*\mu,$ where 
$$u*\mu = \int \limits_{y\in \widehat{G}} u(x-y) d\mu(y)$$
So we just need to show that for some Haar measure $m_H$ on H
$$A=\int \limits_{x\in H} \int \limits_{y\in \widehat{G}} u(x-y) d\mu(y) dm_H(x) < \infty$$

Changing the order of summation in the non-negative integral, 

$$A =\int \limits_{y\in \widehat{G}} \int \limits_{x\in H} u(x-y) dm_H(x) d\mu(y) $$
Let us denote $v=\pi*u^{|^{\widehat{G}/H}},$ where $\pi: \widehat{G} \rightarrow \widehat{G}/H $ is the natural homomorphism. By Theorem 4.1 (3) the above formula implies that
$$A \le const \int \limits_{y\in \widehat{G}} v(-y) d\mu(y)$$
Because $v\in NPPDf(H),$ $v(-y)=v(y) \le 1.$ So
$$A\le const\cdot \int \limits_{y\in \widehat{G}} d\mu (y) = const \cdot f (o) < \infty$$
\hfill \checkbox

\section{Open questions}

There are several open questions that arise naturally from what have proven in the paper.

{\bf Question 1.} Is it true that in Theorem 4.1 (3) the ``almost all" is actually ``all". Or, maybe one can cook up an example of a good function whose average over some coset is strictly less than the limit of the averages over the nearby cosets? My intuition does not tell me what the answer should be, so I am not going to make any conjectures. I just hope that some specialist could provide an answer.

{\bf Question 2.} What is the structure of the cones of good and PPD functions on an arbitrary locally compact abelian group $G$? As far as I know this question is wide open. The following two results suggest that it could be quite interesting.

\begin{Theorem} If $|G| =n < \infty$ then good functions form an open cone in the space $\R^n$ of all real-valued functions on $G.$ The cone of $PPD$ functions on $G$ is the inside the closure of the cone of good functions. This closure is a polyhedral cone.
\end{Theorem}

{\bf Proof.} The cone of good functions is open because it is cut out by a finite number of strict linear inequalities in $R^n.$ The same inequalities cut out the cone of PPD functions, but now those of them that correspond to the elements in $\widehat{G}$ are not strict. \hfill \checkbox

It could be interesting to figure out the combinatorial structure of this cone. By the way, in the natural basis formed by the characteristic functions of the points this cone is algebraic. This means that the coefficients of the defining hyperplanes can be chose to be algebraic numbers. In fact, they can be chosen from the totally real field $\Q [\zeta_n + \zeta_n^{-1}],$ where $\zeta_n$ is a primitive $n-$th root of unity.

\begin{Theorem} Suppose $G=\R^n, $ and $Q$ is a positive-definite quadratic form on $\R^n.$ Then the function $f(x)=e^{-Q(x)}$ is good and generates an extremal ray in the cone of good functions.
\end{Theorem}

{\bf Proof.} If $n=1$ this follows from the theorem of Hardy (cf. \cite{Hardy}). For the higher $n$ it follows from the generalization of the Hardy's theorem, due to Sitaram, Sundari, and Thangavelu (cf. \cite{SST}, also \cite{FS}). I owe the above references to Joaquim Ortega-Cerda. \hfill \checkbox

This theorem of Hardy and its generalizations can be viewed as some particular forms of the general Uncertainty Principle in harmonic analysis. This principle suggests that a function and its Fourier transform could not be simultaneously sharply localized. The above theorem is clearly just the first step, and the structure of the cone of good functions on $\R^n$ is still a mystery. We don't even know if there are any other extremal rays. It is also interesting to find out if this kind of phenomenon happens for any other groups.

{\bf Question 3.} Locally compact abelian groups have many different generalizations, e.g. general locally compact groups, hypergroups, etc. A lot of  harmonic analysis, in particular Pontryagin duality, has been carried out for many such generalizations. So it is natural to ask for the generalizations of the results of this paper to these more general objects. The notion of PPD functions and measures is not hard to make sense of. But it is not clear, e.g. what a right definition of a good function on a non-commutative locally compact group is. This sort of questions is beyond my area of expertise. So I just hope that some harmonic analysis specialists will eventually find the natural framework for the results of this paper.


\begin{thebibliography}{99}

\bibitem{Bor}  Borisov, A. Convolution structures and arithmetic cohomology, preprint,
http://www.math.psu.edu/borisov/researchpage.html . 

\bibitem{Fol} Folland, G. B. A course in abstract harmonic analysis. {\it Studies in Advanced Mathematics} CRC Press, 1995.
 
\bibitem{Hardy} Hardy, G.H. A theorem concerning Fourier transforms, {\it J. London Math. Soc.} {\bf 8} 227--231.

\bibitem{SST} Sitaram, A., Sundari, M., and Thangavelu, S. Uncertainty principles for certain Lie groups, {\it Proc. Indian Acad. Sci. Math. Sci.} {\bf 105} 135 --151. 
\bibitem{FS} Folland, G. B., Sitaram, A. The Uncertainty Principle: a mathematical survey, {\it J. Fourier Anal. Appl.} {\bf 3} (1997), no3, 207--238. 
\end{thebibliography}
\end{document}